\documentclass[a4paper]{amsart}

\usepackage{amssymb,enumerate,amsmath}

\usepackage[utf8x]{inputenc}
\usepackage[english]{babel}
\usepackage[T1]{fontenc}
\usepackage[top=4cm, bottom=3.5cm, left=3.5cm, right=3.5cm]{geometry}

\usepackage{enumerate}
\usepackage[colorlinks]{hyperref}
\usepackage{tikz}
\usepackage{mathtools}
\usepackage[colorlinks]{hyperref}
\usepackage{mathrsfs} 

\newtheorem{thm}{Theorem}[section]
\newtheorem{corollary}[thm]{Corollary}
\newtheorem{prop}[thm]{Proposition}
\newtheorem{lemma}[thm]{Lemma}
\newtheorem{fact}[thm]{Fact}

\newtheorem{conj}{Conjecture}

\theoremstyle{definition}
\newtheorem{defn}[thm]{Definition}
\newtheorem{example}[thm]{Example}

\theoremstyle{remark}
\newtheorem{remark}[thm]{Remark}

\newcommand{\bt}{\begin{thm}}
\newcommand{\et}{\end{thm}}
\newcommand{\bp}{\begin{prop}}
\newcommand{\ep}{\end{prop}}
\newcommand{\bd}{\begin{defn}}
\newcommand{\ed}{\end{defn}}
\newcommand{\bl}{\begin{lemma}}
\newcommand{\el}{\end{lemma}}
\newcommand{\bfa}{\begin{fact}}
\newcommand{\efa}{\end{fact}}
\newcommand{\bc}{\begin{corollary}}
\newcommand{\ec}{\end{corollary}}
\newcommand{\bex}{\begin{example}}
\newcommand{\eex}{\end{example}}
\newcommand{\br}{\begin{remark}}
\newcommand{\er}{\end{remark}}
\newcommand{\ben}{\begin{enumerate}}
\newcommand{\een}{\end{enumerate}}

\newcommand{\isom}{\cong}

\newcommand{\sotto}[2]{#1_{#2}}

\newcommand{\rrr}{\rightarrow}

\newcommand{\ideal}[1]{\sotto {{\mathcal I}}{#1}}

\newcommand{\exact}[3]
{0 \rrr #1 \rrr #2
\rrr #3 \rrr 0}

\newcommand{\pso}{\mathbb{P}^3}

\newcommand{\PP}{\mathbb{P}}

\newcommand{\Z}{\mathbb{Z}}

\newcommand{\coo}{{\mathcal O}}

\newcommand{\caf}{{\mathcal F}}

\newcommand{\lcal}{{\mathcal L}}

\newcommand{\HH}{\mathrm{H}}
\newcommand{\WT}{\textsc{\small WT}}

\newcommand{\Hilb}[2]{\mathbf{Hilb}_{#1,#2}}
\newcommand{\Hlcm}[2]{\mathbf{H}_{#1,#2}}

\begin{document}

\renewcommand{\baselinestretch}{1.2}

\title[The maximum genus problem for locally Cohen-Macaulay space curves]{The maximum genus problem for\\ locally Cohen-Macaulay space curves}

\author[V.~Beorchia]{Valentina Beorchia}
\address{Valentina Beorchia\\ Dipartimento di Matematica e Geoscienze\\ Universit\`a degli Studi di Trieste\\ 
         Via Valerio 12/B \\ 34127 Trieste \\ Italy.}
\email{\href{mailto:beorchia@units.it}{beorchia@units.it}}
\urladdr{\url{http://www.dmi.units.it/~beorchia/}}

\author[P.~Lella]{Paolo Lella}
\address{Paolo Lella\\ Dipartimento di Matematica\\ Politecnico di Milano\\ 
         Piazza Leonardo da Vinci 32\\ 20133 Milano\\ Italy.}
\email{\href{mailto:paolo.lella@polimi.it}{paolo.lella@polimi.it}}
\urladdr{\url{http://www.paololella.it/}}

\author[E.~Schlesinger]{Enrico Schlesinger}
\address{Enrico Schlesinger\\ Dipartimento di Matematica\\ Politecnico di Milano\\ 
         Piazza Leonardo da Vinci 32\\ 20133 Milano\\ Italy.}
\email{\href{mailto:enrico.schlesinger@polimi.it}{enrico.schlesinger@polimi.it}}
\urladdr{\url{http://www1.mate.polimi.it/\%7Eenrsch/}}


\begin{abstract}
Let $P_{\textsc{max}}(d,s)$ denote the maximum
arithmetic genus of a locally Cohen-Macaulay curve  of degree $d$ in $\PP^3$ that is not contained in a surface of degree $<s$. A bound $P(d, s)$ for $P_{\textsc{max}}(d,s)$ has been proven by the first author in characteristic zero and then generalized in any characteristic by the third author. In this paper, we construct a large family $\mathcal{C}$ of primitive multiple lines and we conjecture that the generic element of $\mathcal{C}$ has good cohomological properties. With the aid of \emph{Macaulay2} we checked the validity of the conjecture for $s \leq 100$.  From the conjecture it would follow that  $P(d,s)= P_{\textsc{max}}(d,s)$ for $d=s$ and for every $d \geq 2s-1$.
\end{abstract}



\thanks{This research is supported by MIUR funds PRIN 2015 project {\it Geometria delle variet\`a algebriche} (coordinator A.~Verra) and by MIUR funds FFABR-BEORCHIA-2018. All authors are members of GNSAGA}
\maketitle

\section{Introduction}

\subsection{The maximum genus problem for smooth space curves}
The {\em genus} of a compact connected Riemann surface $X$ is the dimension of
the space $\HH^0 (X, \Omega_X)$ of holomorphic one forms on $X$. It is the unique
topological invariant of $X$.
The set of compact Riemann surfaces of genus $g$ has the structure of an irreducible quasi-projective algebraic variety $\mathcal{M}_g$, which is known as the {\em moduli space}
of surfaces of genus $g$.

A compact Riemann surface can always be embedded in some projective space $\PP^N_{\mathbb{C}}$; its image is then
a {\em smooth projective curve} $C$. The genus of $X$ can be computed by looking at the vector space of hypersurfaces
of degree $n \gg0$ containing $C$: indeed, if one denotes by  $\HH^0 (\PP^N, \coo_{\PP^N}(n))$ the space of all
hypersurfaces in $\PP^N$ of degree $n$ and by $\HH^0 (\PP^N\!, \ideal{C}(n))$ the space of those containing $C$, then
for $n \gg 0$
\begin{equation}\label{hilb-polyn}
\dim \HH^0 \big(\PP^N, \coo_{\PP^N}(n)\big) - \dim \HH^0 \big(\PP^N, \ideal{C}(n)\big)= \deg (C) n + 1- p_a (C).
\end{equation}
The left hand-side of \eqref{hilb-polyn} for $n \gg0$ is thus equal to a degree $1$ polynomial in $n$ that is called the {\em Hilbert polynomial} of $C$;
the coefficient $\deg (C)$ of $n$ is called the {\em degree} of $C$ and it can be computed as the number of intersections of
$C$ with a general hyperplane, while  $p_a (C)$ is called the {\em arithmetic genus} of $C$. It is a theorem that
the arithmetic genus $p_a(C)$ coincides with the topological genus of $X$.

In the algebraic context, over an arbitrary algebraically closed field $k$, one can still define as above the Hilbert polynomial, hence
the degree and arithmetic genus, of an arbitrary, not necessarily smooth, one dimensional subscheme $C \subset \PP^N_{k}$  by
proving that the left hand-side of  \eqref{hilb-polyn} is a degree $1$ polynomial in $n$. Grothendieck has constructed what he called the
{\em Hilbert scheme} $\Hilb{d}{g} (\PP^N_{k})$ of curves of degree $d$ and arithmetic genus $g$ in $\PP^N_k$; the Hilbert scheme $\Hilb{d}{g} (\PP^N_{k})$
parametrizes one dimensional subschemes of degree $d$ and arithmetic genus $g$ in $\PP^N_k$ in the sense that there is a one to one correspondence
between the set of closed points of $\Hilb{d}{g} (\PP^N_{k})$ and the set of such subschemes.

Thus the problem of classifying embedded projective curves became that of determining the pairs $(d,g)$ for which  $\Hilb{d}{g} (\PP^N_{k})$
is non-empty, and then describe $\Hilb{d}{g} (\PP^N_{k})$ for these pairs.
For $N=3$, the problem of determining the possible
pairs $(d,g)$ for smooth connected curves was settled by Gruson and Peskine \cite{gp2}; the correct answer was known already to Halphen  \cite{halphen}
but he lacked a correct proof.

Unfortunately, unlike the moduli spaces $\mathcal{M}_g$, the Hilbert schemes $\Hilb{d}{g} (\PP^N_{k})$ do not have a nice geometric structure except for the fact that they are connected \cite{hthesis};
in fact they may fail to be reduced \cite{extremal,mumford-fpathologies}  or irreducible, and may have all type of singularities \cite{Vakil-Murphys-Law}. 
Therefore any serious study of embedded curves involves more refined invariants than just degree and genus to obtain a stratification
of $\Hilb{d}{g} (\PP^N_{k})$ in nicer subschemes: this was already clear to Halphen \cite{halphen}. A natural choice is to
consider the {\em postulation} of $C$, that is, the dimension of the space of hypersurfaces of degree $n$ containing $C$, for every $n$.
By \eqref{hilb-polyn} the postulation refines the degree and genus. As it is very difficult to
characterize the possible values of the postulation for curves of a given degree and genus, beginning with Halphen
researchers have concentrated their efforts on the intermediate problem of determining the {\em maximum genus $G(d,s)$ of a curve of degree $d$ in $\PP^3$ that does not lie
on a surface of degree $<s$}. This was certainly motivated by the observation that curves of very high genus have to be contained in surfaces of very low degree.
The problem is now known as the {\em maximum genus problem}, and it is not completely settled yet: the pairs $(d,s)$ to be considered are naturally
divided in three domains called range $A$, $B$ and $C$ respectively \cite{MR89g:14025}. In range $C$, the precise value of $G(d,s)$ has been determined  by Gruson and Peskine
\cite{gp1}.  In the most difficult range $B$,
there is a conjectural bound $G_B(d,s)$ for $G(d,s)$ and Hartshorne and Hirshowitz \cite{hhnouvelles} have constructed curves of degree $d$ and genus $G_B(d,s)$ not lying on surfaces of degree $<s$. It remains to show that $G(d,s) \leq G_B(d,s)$. There is recent very interesting work of Macri and Schmidt on this topic \cite{macri-schmidt}.

To give the reader the flavour of the maximum genus problem, we would like to explain the conjectural value of $G(d,s)$ in range $A$, namely when  $d < \frac{1}{3}  (s^2+4s+6)$. Let
\[
G_A (d,s)=  (s-1) d + 1 - \dim \HH^0 \big(\PP^3, \coo_{\PP^3}(s{-}1)\big) = (s-1) d + 1 - \binom{s+2}{3}.
\]
The inequality $d < \frac{1}{3}  (s^2+4s+6)$ together with Clifford's theorem implies $\HH^1 (C, \coo_C(s{-}1))=0$; therefore in range $A$
a curve of degree $d$ not lying on a surface of degree $s-1$ must have genus $g(C) \leq G_A(d,s)$,  and equality holds
if and only if $C$ has maximal rank --- recall that a
curve $C$ in $\PP^N$ is said to be of maximal rank if the natural map $\rho(n): \HH^0 (\PP^N, \coo_{\PP^N}(n)) \rightarrow \HH^0 (C, \coo_{C}(n))$ has maximal rank, that is,
is either surjective or injective, for every integer $n$. Thus, if one could construct  for each $(d,s)$ satisfying
$d < \frac{1}{3}  (s^2+4s+6)$  an irreducible smooth curve of maximal rank, degree $d$ and genus $G_A(d,s)$,
it would follow that $G_A(d,s)=G(d,s)$.  Such curves have been constructed in a subinterval of range A (see \cite{BBEM,BEF}). The curves constructed by Gruson-Peskine to prove sharpness in range C and the examples by Hartshorne-Hirshowitz in range B are curves of maximal rank as well.


We would like to mention here that recently Eric Larson has proven the
maximal rank conjecture for a Brill-Noether general curve in a rather extraordinary series of papers culminating in \cite{larson} --- the case
$N=3$ of the conjecture was settled by Ballico and Ellia in the 80's \cite{ballico-ellia}.

\subsection{The maximum genus problem for locally Cohen-Macaulay space curves}
If on the one hand the class of {\em smooth curves} is too restricted for most purposes, on the other hand the class
of one dimensional subschemes with a fixed Hilbert polynomial is too large if one is interested in the geometry of curves.
From the point of view of liaison theory, the correct class to be considered it that of  locally Cohen-Macaulay
curves, that is, one dimensional subschemes of $\PP^3$ with no zero-dimensional irreducible or embedded components.
Locally Cohen-Macaulay curves form an open subset $\Hlcm{d}{g}$ of the full Hilbert scheme $\Hilb{d}{g} (\PP^3) $, and
one can use liaison theory and the Hartshorne-Rao module $\bigoplus_{n \in \Z} \HH^1 (\pso, \ideal{C}(n))$ to obtain a nice stratification
of $\Hlcm{d}{g}$: the set of locally Cohen-Macaulay space curves of degree $d$ and genus $g$ curves with fixed cohomology and Hartshorne-Rao module is \emph{smooth} and \emph{irreducible} and its dimension is in principle computable
\cite[Th\'eor\`{e}me VII1.1, Corollaire VII.1.7, Proposition 3.2]{MDP}.

The problem of determining the possible values $(d,g)$ for the degree and arithmetic genus of a locally Cohen-Macaulay space curve
turned out to be much easier than the corresponding problem for smooth curves. Indeed, the Hilbert schemes are nonempty
if and only if $d \geq 1$ and $g=(d-1)(d-2)/2$ or if $d \geq 2$ and $g \leq (d-2)(d-3)/2$ (see \cite{hgenus}).

This paper is concerned with the problem of maximum genus for locally Cohen-Macaulay space curves: determine the maximum
arithmetic genus of a locally Cohen-Macaulay space curve of degree $d$ that is not contained in a surface of degree $s-1$.
The first author \cite{beo2} proved a bound $P(d,s)$ for the maximum genus if the characteristic of the ground field is zero, and proved the bound is sharp if
$s \leq 4$; later the third author \cite{sch-beo} gave a different proof of this bound valid in any characteristic. Note that there is
a locally Cohen-Macaulay curve of degree $d$ that is not contained in a surface of degree $<s$ if and only if
$ d \geq s \geq 1$. It is clear that one must have $d \geq s$ because a curve is contained in the surface obtained as cone over
a general plane section, while, if $d \geq s$, an example of a curve of degree $d$ not contained in a surface of degree
$<s$ is the divisor $C=dL$ on $S$ where $L$ is a line contained in a smooth surface $S$ of degree $s$.

\bt[\cite{beo2,sch-beo}]
\label{beo}
Let $C$ be a curve in $\pso$ of degree $d$ and genus $g$. Assume
that $C$ is not contained in any surface of degree $<s$.
Then $d \geq s$ and
\begin{equation}\label{Bbound}\tag{$\star$}
g \leq P(d,s) =
\begin{cases}
(s-1)d +1 - \binom{s+2}{3}, & \text{if~}  s \leq d \leq 2s, \\
\binom{d-s}{2}-
\binom{s-1}{3},& \text{if~} d \geq 2s+1.\\
\end{cases}
\end{equation}
\et

Thus to show that $P(d,s)$ is the maximum
arithmetic genus $P_{\textsc{\textsc{max}}}(d,s)$ of a locally Cohen-Macaulay space curve of degree $d$ that is not contained in a surface of degree $s{-}1$
it remains to construct, for each pair $d \geq s$, a curve of genus $P(d,s)$ not lying on a surface of degree $s{-}1$.
Note that numerically $P(d,s)=G_A(d,s)$, but the domains of these functions in the smooth irreducible case and in the locally Cohen-Macaulay case are different.

We think the case $d=s$ is the crucial case of the problem. Indeed,  if $P(s{-}1, s{-}1)=P_{\textsc{max}} (s{-}1,s{-}1)$,
it easily follows that $P(d,s)=P_{\textsc{max}}(d,s)$ for every
$d \geq 2s-1$ (see Section \ref{section 5}). So let us suppose $d=s$. It is known that curves of degree $d$ not lying on a surface of degree $d-1$
are supported on one line or on two disjoint lines \cite{sch-tran}. Among curves supported on a single line, the ones  with the nicest scheme structure are those contained locally in a smooth surface; they are called {\em primitive multiple lines} \cite{banica,drezet}. Let $C$ be a primitive multiple line of degree $d$ having the line $L$ as its support. Then the sheaf
$\lcal=\ideal{L,C}/\ideal{L,C}^2$ is an invertible sheaf of $\coo_L$-modules; indeed,  over the open set where $C=dL$ on a smooth surface $S$,  $\lcal$ is the conormal sheaf of $L$ in $S$. Therefore there exists a unique integer $e$ such that $\lcal \cong \coo_L(e)$. We call
$e$ the type of the primitive $d$-line $C$. 
For $e \geq 0$, by a result of Beorchia and Franco \cite[Proposition 1.1]{BF}, such a $C$ is an embedding in $\mathbb{P}^3$ of $dC_0 \subseteq \mathbb{F}(e)$, where $\mathbb{F}(e)$ is the ruled surface $\mathbb{P}_{\mathbb{P}^1} (\mathcal{O}_{\mathbb{P}^1}\oplus\mathcal{O}_{\mathbb{P}_1} (e))$ and $C_0$ is a section with $C_0^2 = - e$; see also \cite{drezet}.

One easily computes that the genus of a primitive $d$-line $C$ of type $e$ is
\[
p_a (C)=-\sum_{i=1}^{d-1}(ie+1).
\]
In order to show the bound $P(d,d)$ is sharp,  it is then enough to prove the following conjecture about existence of primitive $d$-lines with so to speak {\em good cohomology}, an analogue of the maximal rank condition.
 Note that a primitive $d$-structure $C$ on a line $L \subset \pso$ is contained in the
infinitesimal neighborhood $L_d$ of $L$ defined by the $d$-th power $\ideal{L,\pso}^d$ of the ideal sheaf of
$L$.
%
%


\begin{conj} \label{A}
Let $d \geq 5$ be an integer congruent to $2$ modulo 3. Then
\ben[(i)]
\item there is a primitive $d$-line $C$ of type $e=\frac{d-2}{3}$ that is not contained in a surface
of degree $d-1$;
\item
there is a primitive $d$-line $C$ of type $e=\frac{d+1}{3}$ with support $L$ such that
\[
\HH^0 \big(\ideal{C}(n)\big)= \HH^0 \big(\ideal{L}^d (n)\big)\quad \text{for~} n \leq d.
\]
This means that $C$ is not contained in a surface of degree $d-1$ and furthermore, every surface of degree
$d$ containing $C$ contains $L_d$;
\item
there is a primitive $d$-line $C$ of type $e=\frac{d+4}{3}$ with support $L$  such that
\[
\HH^0 \big(\ideal{C}(n)\big)= \HH^0 \big(\ideal{L}^d (n)\big) \quad \text{for~} n \leq d+1.
\]
\een
\end{conj}

Conjecture \ref{A} implies that the bound $P(d,d)$ for the genus of a curve of degree $d$ not lying on surfaces of degree $<d$ is sharp.
Indeed, if $d \equiv 2 \pmod 3 $, a primitive $d$-line $C$ of type $e=\frac{d-2}{3}$ has genus $P(d,d)$ and, according to Conjecture \ref{A}
there is such a curve $C$ not lying on surfaces of degree $<s$.

If
$d \equiv 0 \pmod 3$, the disjoint union $D$ of a line $L'$ and a primitive $(d-1)$-line $C$ of type $e=\frac{d}{3}$ has
genus $P(d,d)$; by Conjecture \ref{A} there is such a curve $C$ with support $L$ satisfying
$\HH^0 (\ideal{C}(d-1))= \HH^0 (\ideal{L}^{d-1} (d-1))$; so every surface of degree $<d$ containing $C$ is a union of planes containing $L$. Thus,
if we choose a line $L'$ disjoint from $L$, no surface of degree $d$ can contain $L' \cup C$.

Finally, if $d \equiv 1 \pmod 3$, $e=\frac{d+2}{3}$, the disjoint union $D$ of a double line $B$ of type $d\!-\!2$ and a primitive $(d\!-\!2)$-line of type
$e=\frac{d+2}{3}$ has genus $P(d,d)$. By Conjecture \ref{A} there is such a curve $C$ with support $L$ satisfying
$\HH^0 (\ideal{C}(d-1))= \HH^0 (\ideal{L}^{d-2} (d-1))$. On the other hand,
if the support of the double line $B$ is the line $L'$ of equations $x=y=0$, the homogenous ideal of $B$ is generated by $x^2,xy,y^2$ and by the equation of
a surface of degree $d=e(B)+2$. The condition $\HH^0 (\ideal{C}(d-1))= \HH^0 (\ideal{L}^{d-2} (d-1))$ implies there is no nonzero  $F \in \HH^0 (\ideal{C}(d-1))$
that is contained in the ideal $(x,y)^2$ if the lines $L'$ and $L$ are disjoint. Thus no surface of degree $d$ can contain the curve $D=B \cup C$ if the supports
of $B$ and $C$ are disjoint.

The main purpose of this paper is to introduce a construction of primitive multiple lines that we hope will lead to a proof of the conjecture.
We begin by motivating our construction. As
a primitive $d$-structure $C$ on a line $L \subset \pso$ is contained in the
infinitesimal neighborhood $L_d$ of $L$, we can consider the ideal sheaf $\caf(C)$ of $C$ in $L_d$. Since $C$ is primitive, the sheaf
$\caf(C)$ turns out to be an invertible sheaf of $\coo_{L_{d-1}}$-modules: a local generator of $\caf(C)$ at a point $p$
is the image of the local equation $g \in \ideal{C,p}$ of a smooth surface containing $C$ near $p$.
The Hilbert polynomial of this sheaf is $
\chi {\caf}(n)=   \left(n-e-\frac{1}{3}(2d-1)\right)\binom{d}{2}$. Hence the condition $d \equiv 2 \pmod 3$
guarantees that the unique zero $n_0$ of $\chi {\caf}(n)$ is an integer, so that to prove the conjecture it is enough to construct
 a primitive $d$-line $C$ of type $e$ such that $\HH^0 (\caf(n_0))=0$. Note that $\HH^0 (\caf(n_0))=0$ implies that $\caf$ has natural cohomology,
 that is, $ \dim \HH^0 (\caf(n)) = \max (\chi {\caf}(n),0)$ for all $n \in \Z$. 

Furthermore, if $S \subset \pso$ is a surface containing $L$ which is generically smooth along $L$,
then $S$ contains a unique locally Cohen-Macaulay $d$-structure $C$ on the line $L$ because, being Cohen-Macaulay,
such a curve is determined by its local ring at the generic point of $L$, where $C$ is the intersection of $S$
and $L_d$. The equation of $S$ defines a section of the invertible sheaf $\caf(C)$, hence a divisor $E$ on $L_{d-1}$.
The condition $\HH^0 (\caf(n_0))=0$ translates on a condition on the cohomology of $E$.

An analysis of the possible local equations for such an $E$ led us to consider polynomials in three variables $x$, $y$, $z$
that are homogeneous with respect to the weight function \WT\ that assigns weight $1$ to $x$, weight $2$ to $y$ and  weight $3$ to $z$.
Now fix an integer $d=3m-1 \geq 5$ congruent to $2$ modulo $3$. A \WT-homogeneous polynomial of weight $3m$ has the form
\[
g(x,y,z)=z^m + \sum_{i=1}^m a_i(x,y)z^{m-i}.
\]
We assume that $g$ has no term in $(x,y)^{d-1}$. Every monomial in $x,y,z$ of weight $3m$ that does not belong to $(x,y)^{d-1}$ has degree strictly less than $d-1$,
except for $x^{3(m-1)}z$ that has degree $d-1$. We assume that   $x^{3(m-1)}z$ appears in $g$ with nonzero coefficient, so that $g$ has degree $d-1$.
Let $G(X,Y,Z,W)= W^{d-1} g(X/W,Y/W,Z/W)$ be the homogenization of $g$ with respect to the standard grading; it is a polynomial of degree $d-1$. The main result of this paper is the following:

\bt \label{main result}
In the projective space $\pso = \textnormal{Proj} (k[X,Y,Z,W])$ over the algebraically closed field $k$, consider the line $L$ of equations $X=Y=0$ and, for an integer $a \geq 0$ and a scalar $t \in k^*$,
the surface $S_t$ of equation $XGZ^a-tYW^{a+d-1}=0$. Let $C_t$ be the unique locally Cohen-Macaulay $d$-structure on $L$ contained in the surface $S_t$.
Then
\begin{enumerate}[ (i)]
\item\label{main:a} for every integer $a \geq 0$ and a scalar $t \in k^*$, the curve
$C_t$ is a primitive $d$-line of type $e=a+m-1$;
\item\label{main:b}
suppose furthermore that the ideal $I=(x,y)^{d-1}+(g) \subset k[x,y,z]$ contains no nonzero polynomial of degree $d-2$. Then
\begin{equation*}\label{cohomology}
\HH^0 (\ideal{C_t}(n))= \HH^{0} (\ideal{L}^d(n)) \qquad \mbox{for all $n \leq a+d-1$}.
\end{equation*}
\end{enumerate}
In particular, if a polynomial $g$ as in (\ref{main:b}) exists,  Conjecture \ref{A} is true for $d=3m-1$, and the bound $P(d',d')$ is sharp for $d'=d$, $d+1$ and $d+2$: there exists
a curve of degree $d'$ and arithmetic genus $P(d',d')$ that does not lie on a surface of degree $<d'$.
\et

In light of Theorem \ref{main result}, Conjecture \ref{A} is thus a consequence of the following

\begin{conj}\label{conjB}
Fix an integer $m \geq 2$ and let $g$ be a general \WT-homogeneous polynomial in $k[x,y,z]$ of weight $3m$. Then the ideal $I=(x,y)^{3m-2}+(g)$ contains no polynomial of standard degree $3m-3$.
\end{conj}

We have verified Conjecture \ref{conjB} for $m=2$ and $m=3$ by handwritten computations and with the aid of \emph{Macaulay2} \cite{M2} for $m \leq 40$ and  a ground field $k$ of characteristic sufficiently large.


In Section \ref{section gpol} we prove that the $k$-algebras $k[x,y,z]/(x,y,z)^{d-1}$ and $k[x,y,z]/((x,y)^{d-1}$ ${}+(z^m))$ have the same Hilbert function with respect
to the \WT-grading; if the Hilbert scheme parametrizing such zero dimensional (but weighted) affine schemes with fixed \WT-Hilbert function were irreducible,
the generic point would give a deformation $g$ of $z^m$ as desired. Or one would be done if one could show that for general
$g$ of weight $3m$ the initial ideal of  $(x,y)^{3m-2}+(g)$ is $(x,y,z)^{3m-2}$, with respect to some (or any) monomial order
that refines the usual degree. 

%

\subsection{Content of the paper}
In Section  \ref{section 2} we review some basics fact about primitive multiple lines in $\pso$. In Section  \ref{section 3}, Proposition \ref{liftDtoC},
we prove that the curves appearing in Theorem \ref{main result} are indeed primitive multiple lines of the required type; for this we need the condition that $g$ is \WT-homogeneous.
In Section  \ref{section 4} we prove part \emph{(\ref{main:b})} of Theorem  \ref{main result} reducing the statement about the cohomology of $\mathcal{F}(C) = \mathcal{I}_C/\mathcal{I}_{L_d}$ to a statement about the zero dimensional scheme defined
by the ideal $I=(x,y)^{d-1}+(g)$. In Section \ref{section gpol} we study what this condition means in terms of the polynomial $g$ and we explain how we checked Conjecture \ref{conjB} for $m \leq 40$.
In the final Section \ref{section 5} we prove that, if $P(s{-}1,s{-}1)$ is sharp, then $P(d,s)$ is sharp for every $d \geq 2s-1$,
thus proving sharpness of bounds \eqref{Bbound} in most cases assuming Conjecture \ref{conjB} is true.

\subsection*{Acknowledgements} The authors would like to thank Lea Terracini and Michele Rossi for fruitful discussions.

\section{Primitive multiple lines embedded in $\pso$} \label{section 2}
We would like to guess what the postulation of the {\em general} primitive $d$-line in
$\pso$ should be. A restraint on the possible values of the postulation comes from the fact
that any primitive $d$-line $C$ in $\pso$ is contained in the $(d-1)$-th infinitesimal neighborhood
$L_d$ of $L=C_{\textnormal{red}}$, that is,
\[
\ideal{L}^d \subseteq \ideal{C}.
\]
Thus it is natural to look at the sheaf
\[
\caf=\caf(C):=\ideal{C,L_d}=\ideal{C}/\ideal{L}^d
\]
For an invertible sheaf $\caf$ of $\coo_{L_{d-1}}$-modules, we define the {\em reduced degree} of $\caf$ as
the degree of the restriction $\caf_{L}$ of $\caf$ to the reduced line $L$.

\bp \label{basic}
A subscheme $C \subset L_d$ is a primitive $d$-line if and only if $C$  has no zero dimensional
embedded or isolated components, the sheaf $\caf(C)=\ideal{C,L_d}$ is annihilated by
$\ideal{L}^{d-1}$ and $\caf(C)$ is locally free of rank $1$, as a module over $\coo_{L_{d-1}}$.
If $C$ is primitive and $e$ is its type, the reduced degree of $\caf(C)$ is $-e-2$. At a point $p$,
if $f \in \ideal{C,p}$ defines locally a smooth surface containing $C$, then the image of $f$ in  $\caf(C)_p$ generates $\caf(C)_p$.
\ep

\begin{proof}
The sheaf $\caf(C)=\ideal{C,L_d}$ is contained in $\ideal{L,L_d}$, thus it is
annihilated by $\ideal{L}^{d-1}$. Locally, $C$ is contained in a smooth surface $S$,
so the ideal of $C$ is locally generated by a regular sequence $x$, $y^d$ where $x$ is the equation
of $S$ and $x$, $y$ are local equations of $L$. Then $[x]$ is a free local generator of $\caf(C)=\ideal{C,L_d}$
over $\coo_{L_{d-1}}$ because the annihilator of $[x]$ is $\ideal{L}^{d-1}$.
Finally, to compute the type of $C$, we consider the unique $2$ structure $C_2$ on $L$ contained in $C$: if
$f \in \ideal{C,p}$ defines locally a smooth surface containing $C$, then the ideal of $C_2$ at $p$ is generated
by $f$ and $\ideal{L}^2$. It follows that  \[\caf (C) \otimes \coo_L \cong \caf (C_2) \otimes \coo_L \cong \caf(C_2).\]
On the other hand, the type $e$ of $C$ is the same as the type of $C_2$, and for $C_2$ we use the exact sequence
\[
\exact{\caf(C_2)=\ideal{C_2}/\ideal{L}^2}{\ideal{L}/\ideal{L}^2}{\ideal{L}/\ideal{C_2}\isom \coo_L(e)}
\]
to deduce $\caf(C_2) \cong \coo_L(-e-2)$.
\end{proof}

\bex
If $C$ is contained in a surface $S$ of degree $s$ smooth along $L=C_{\textnormal{red}}$, then
\[
\caf(C) \cong \coo_{L_{d-1}}(-s)
\]
and  the type of $C$ is $e=s-2$. One can compute the type of $C$ most easily by finding a smooth
surface $S$ containing $C_2$, then the type of $C$ is $\deg(S)-2$.
\eex

More generally:
\bl \label{lemma}
Suppose $C \subset \pso$ is a primitive $d$-line of type $e$ with support $L$, and suppose
$S$ is a surface of degree $s$ that contains $C$ and is generically smooth along $L$. Then the equation of $S$ gives rise to
an exact sequence
\[
\exact{\coo_{L_{d-1}}(-s)}{\ideal{C}/\ideal{L}^d}{\coo_E}
\]
where $E$ is an effective Cartier divisor on $L_{d-1}$ of reduced degree $s-e-2$ that satisfies
$\coo_{L_{d-1}}(E) \cong \ideal{C}/\ideal{L}^d (s)$.
\el

\begin{proof}
Since $S$ is generically smooth along $L$, its equation gives a regular section of the invertible $\coo_{L_{d-1}}$-module
$\caf (s) = \ideal{C}/\ideal{L}^d(s)$; this section defines the Cartier divisor $E$, and gives rise to the standard exact sequence
\[
\exact{\coo_{L_{d-1}}}{\caf(s)\cong \coo_{L_{d-1}}(E)}{\caf(s) \otimes \coo_E \isom \coo_E}
\]
where the last isomorphism is due to the fact that $E$ is zero dimensional. To compute the reduced degree of $E$ we recall that
the reduced degree of $\caf = \ideal{C}/\ideal{L}^d$ is $-e-2$, hence
the reduced degree of $\caf(s)\cong \coo_{L_{d-1}}(E)$ is $s-e-2$.
\end{proof}

\br
The postulation of $C$ is determined by the $\mathrm{h}^0$-function $n \mapsto \mathrm{h}^0 (\caf(n)) = \dim \HH^0(\caf(n))$ of $\caf$
because $L_d$ is arithmetically Cohen-Macaulay. 
\er

\bp \label{hp}
The Hilbert polynomial of $\caf$ is
\begin{equation*}
\chi {\caf}(n)=        \left(n-e-\frac{1}{3}(2d-1)\right)\binom{d}{2}.
\end{equation*}
\ep
\begin{proof}
When $C$ is contained in a degree $s$ surface smooth along $L=C_{\textnormal{red}}$, then $s = e+2$ and $\caf(C) \cong \coo_{L_{d-1}}(-e-2)$;
since the Hilbert polynomial depends only on $d$ and $e$, we conclude that
$\chi \caf(n)= \chi  \coo_{L_{d-1}}(n-e-2)$, for every primitive $d$-line of type $e$.
\end{proof}

\section{Construction of primitive multiple lines} \label{section 3}
Define a weight function \WT\ on $k[x,y,z]$ setting as above $\WT(x)=1$, $\WT(y)=2$, $\WT(z)=3$. We fix an integer $d \geq 5$ congruent to
$2$ modulo $3$, and we write $d=3m-1$ with $m\geq 2$.  Let
\[
g_0(x,y,z)= z^m+a_1(x,y)z^{m-1}+\cdots + a_{m-1} (x,y)z + a_m(x,y)
\]
denote a \WT-homogeneous polynomial in $k[x,y,z]$ of weight $3m$, and assume all monomials appearing in the coefficients $a_j(x,y)$ have standard degree at most $3(m-1)$;
essentially, we are considering a \WT-homogeneous polynomial of weight $3m$ modulo $(x,y)^{d-1}=(x,y)^{3m-2}$. Note that, once we have excluded monomials
in $(x,y)^{3m-2}$, every other monomial of weight $3m$ has
degree at most $d-1=3m-2$, and the only monomial not in $(x,y)^{3m-2}$ of weight $3m$ and maximum degree $3m-2$ is $x^{3(m-1)}z$. We homogenize $g_0$
with respect to the standard grading: more precisely, we let
\[
G(X,Y,Z,W)= W^{d-1} g_0(X/W,Y/W,Z/W)= Z^{m}W^{2(m-1)} + \cdots
\]
Finally, we dehomogenize $G$ with respect to $Z$, that is, we define
\[
g_{\infty}(x,y,w)= G(x,y,1,w)=w^{2(m-1)}+p_1(x,y) w^{2m-3} + \cdots + p_{2(m-1)}(x,y).
\]

Define a grading $\WT_{\infty}$ on $k[x,y,w]$ setting $\WT_{\infty}(x)=2$, $\WT_{\infty}(y)=1$ and $\WT_{\infty}(w)=3$. Then
$g_{\infty}$ is $\WT_{\infty}$-homogeneous of weight $6(m-1)$ and the coefficients $p_j(x,y)$ are $\WT_{\infty}$-homogeneous of weight $3j$: this is
because every monomial $x^ay^bz^c$ of weight $\WT=a+2b+3c$ in $g_0$ gives rise to a monomial $x^ay^bw^{3m-2-(a+b+c)}$ of weight $\WT_{\infty}=9m-6-(a+2b+3c)$
in $g_{\infty}$.

\bl \label{lemmaprim}
Let $g_0$, $G$ and $g_\infty$ be as above.
Then for any integer $a \geq 0$ there exist polynomials $g_1,h,h_1$ in $k[x,y,w]$ such that for all $t \in k$
\begin{equation}\label{factor}
(g_{\infty}+th)(x-ty(g_1+th_1)) \equiv xg_{\infty}-tyw^{a+d-1} \mod (x,y)^{d}.
\end{equation}
\el

\begin{proof}
We will repeatedly use the following preliminary remark:  a monomial in the first two variables $x$ and $y$
 of $\WT_{\infty}$ greater or equal than $6m-5$ is in the ideal $(x,y)^{d-1}$, as the monomial in $a(x,y)$ of maximum weight not
 in $(x,y)^{d-1}$  is $x^{d-2}=x^{3m-3}$ which has weight $6m-6$.

We begin by constructing a $\WT_{\infty}$-homogeneous polynomial $\overline{g}$ of weight $3m$ such that
$g_\infty \overline{g} \equiv w^{d-1}$ modulo ${(x,y^{d-1})}$. We look for such a $\overline{g}$ among polynomials of the form
\[
\overline{g}= w^{m} + q_1(x,y)w^{m-1}+ \cdots +q_m(x,y)
\]
where the coefficients $q_j(x,y)$ are $\WT_{\infty}$-homogeneous of weight $3j$. Now, set $p_0=q_0=1$,
\begin{equation}\label{gg1}
g_\infty \overline{g} = w^{3m-2} + \sum_{\ell=1}^{m} w^{3m-2-\ell} \left(  \sum_{j=0}^{\ell}  p_{\ell-j} q_j \right)
+\sum_{\ell=1}^{2(m-1)} w^{2(m-1)-\ell} r_{\ell}(x,y).
\end{equation}
Note that $\WT_{\infty}(r_\ell)=3m+3\ell$; so, if a monomial in $r_{\ell}(x,y)$ is not divisible by $x$, then it must be
$y^{3m+3\ell}$, and  $3m+3\ell \geq 3m+3=d+4$. Thus, if we define inductively $q_\ell= -\sum_{j=0}^{\ell-1} p_{\ell-j} q_j$ for $\ell=1,2,\ldots, m$
it follows from \eqref{gg1} that 
\[
g_{\infty} \overline{g} - w^{3m-2} \equiv x k(x,y,w) \mod (y^{d+4})
\]
where $k(x,y,w)$ is a $\WT_{\infty}$-homogeneous of weight $9m-8$ of the form
\[
k(x,y,w) = \sum_{j=0}^{2m-3}c_j(x,y) w^{2m-3-j}.
\]
We compute $\WT_{\infty} (c_j(x,y)) = (9m-8)-(6m-9-3j)=3m+3j+1$. Hence, for $j \geq m-2$, the $\WT_{\infty}$-homogeneous
polynomial $c_j(x,y)$ has weight $\geq 6m-5$, so by our initial remark it belongs to $(x,y)^{d-1}$. We conclude
\[
k(x,y,w) \equiv  \sum_{j=0}^{m-3}c_j(x,y) w^{2m-3-j} \mod (x,y)^{d-1}.
\]

Next we look at the polynomial $f=-\overline{g} k$: it is
$\WT_{\infty}$-homogeneous of weight $6m-2$, and the highest exponent of $w$ appearing in $f$
is at most $3m-3$. Thus we can write  $f$ in the form
\[
f(x,y,w) = \sum_{j=0}^{3(m-1)}e_j(x,y) w^{3m-3-j}  \equiv \sum_{j=0}^{m-3}e_j(x,y) w^{3m-3-j} \mod(x,y)^{d-1}.
\]
To derive the last congruence we used the fact that the coefficients $e_j(x,y)$ have weight  $3m+3j+1$
and thus belong to $(x,y)^{d-1}$ for $j \geq m-2$.

We claim that any $\WT_{\infty}$-homogeneous polynomial $f \equiv  \sum_{j=0}^{m-3}e_j(x,y) w^{3m-3-j}$
can be divided by $g_{\infty}$ modulo $(x,y)^{d-1}$: more precisely, we next construct a $\WT_{\infty}$-homogeneous polynomial $\overline{k}$
of weight $6m-2$ such that $g_{\infty}  \overline{k}=f$  modulo $(x,y)^{d-1}$. Note that the highest exponent of $w$ appearing in
$\overline{k}$ will be at most $m-1$, so $\overline{k}$ can be expanded as $\overline{k} = \sum_{i=0}^{m-1} b_i(x,y) w^{m-1-i}$ where the coefficient
$b_i$ has weight $3m+3i+1$;  hence $b_i$ belongs to $(x,y)^{d-1}$ for $i \geq m-2$, and $\overline{k} \equiv \sum_{i=0}^{m-3} b_i w^{m-1-i}$
modulo  $(x,y)^{d-1}$. Thus,
\begin{equation*}
\begin{split}
g_{\infty}  \overline{k} & {} \equiv \left( w^{2m-2}+ \sum_{j=1}^{2(m-1)} p_j w^{2m-2-j}\right)\left( \sum_{i=0}^{m-3} b_i w^{m-1-i} \right) \equiv {}\\
   & {} \equiv b_0 w^{3m-3} + (b_1+b_0 p_1) w^{3m-4}+\cdots +(b_{m-3}+\cdots+ b_0 p_{m-3} ) w^{2m} \mod (x,y)^{d-1}
   \end{split}
\end{equation*}
and we can choose the $m-3$ coefficients $b_i$ so that $g_{\infty}  \overline{k}=f$  modulo $(x,y)^{d-1}$ for any $f= \sum_{j=0}^{m-3}e_j(x,y) w^{3m-3-j}$.

Finally, we observe that $k  \overline{k} \equiv 0$ modulo $(x,y)^{d-1}$. This is because
$k(x,y,w)$ is congruent to $\sum_{j=0}^{m-3}c_j(x,y) w^{2m-3-j}$ with $\WT_{\infty}(c_j)= 3m+3j+1$, and similarly
$\overline{k} \equiv \sum_{i=0}^{m-1} b_i w^{m-1-i}$ with
$\WT_{\infty}(b_i)= 3m+3i+1$. Thus $\WT_{\infty}(b_ic_j) \geq 6m+2$ for every $i,j$, hence
$k  \overline{k} \equiv 0$ modulo $(x,y)^{d-1}$.

Summing up we have constructed polynomials $\overline{g},k,\overline{k}$ such that the following congruences hold:
\[
  g_{\infty} \overline{g} -w^{d-1}-xk \equiv g_{\infty} \overline{k} +  \overline{g} k \equiv  k  \overline{k} \equiv 0 \mod (x,y)^{d-1}.
\]
To prove the statement now it is enough to set $g_1=w^a \overline{g}$, $h=y w^a k$ and $h_1=yw^a  \overline{k}$ in equation \eqref{factor}.
\end{proof}

\bp \label{liftDtoC}
Let $G$ and $g_{\infty}$ be as above. For every integer $a \geq 0$ and scalar $t \in k^*$, consider\hfill the\hfill line\hfill $L$\hfill of\hfill equations\hfill $X=Y=0$\hfill inside\hfill the\hfill surface\hfill $S$\hfill of\hfill equation\\
$XGZ^a-tYW^{a+d-1}=0$. Let $C$ be the unique locally Cohen-Macaulay $d$-structure on $L$ contained in the surface $S$.
Then $C$ is primitive of type $e=a+m-1$. Furthermore, the divisor $E$ associated to the pair $(C,S)$ as in Lemma \ref{lemma}
is supported at $p_{\infty}=[0:0:1:0]$ and has local equation $g_{\infty}+th$ for some polynomial $h \in k[x,y,w]$ divisible by $y$.
\ep

\begin{proof}
Note that the surface $S$ of equation $XGZ^a-tYW^{a+d-1}=0$ is smooth along $L$ except at the point $p_{\infty}=[0:0:1:0]$,
hence $C$ is primitive except perhaps at $p_{\infty}$.
To see what happens at $p_{\infty}$, we work in the affine open subset $\mathbb{A}^3$ where $Z\neq 0$ with coordinates $x=X/Z$, $y=Y/Z$ and $w=W/Z$.
By Lemma \ref{lemmaprim} the equation $f(x,y,w)$ of $S$ in  $\mathbb{A}^3$ factors as
\[
f \equiv (g_{\infty}+th)(x-ty(g_1+th_1)) \mod (x,y)^{d}.
\]
Now recall that $C$ is constructed by intersecting $S$ with the infinitesimal neighborhood $L_d$ of $L$
and the removing embedded points. As $g_{\infty}+th \equiv w^{2(m-1)}$ modulo $(x,y)$, the polynomial
$g_{\infty}+th$ maps to a nonzero divisor in $\coo_{L_d}$, hence $x-ty(g_1+th_1)$ is contained in the ideal
of $C$ at $p_{\infty}$. As  the surface $x-ty(g_1+th_1)=0$ is smooth at $p_{\infty}$, we conclude that $C$ is primitive,
that the image of $x-ty(g_1+th_1)$ in $\ideal{C,L_d}$ is a local generator, and that the divisor $E$ has local equation
$g_{\infty}+th$ at $p_{\infty}$. Furthermore, since $S$ is smooth along $L$ except at $p_{\infty}$, the divisor $E$ is supported at $p_{\infty}$.

To compute the type of $C$, we observe that the multiplicity $2$ substructure $C_2$ of $C$ is contained in the surface of equation
$XZ^{a+m}{W}^{2m-2}-tYW^{a+3m-2}=0$ because $d-1=3m-2$ and $G \equiv Z^{m}W^{2m-2}$ modulo $(X,Y)^2$. Hence $C_2$ is contained in the
surface of equation $XZ^{a+m}-tYW^{a+m}=0$; as this surface is smooth along $L$ and has degree $a+m+1$, the type of $C$ is
$a+m-1$.
\end{proof}

\section{Reduction to the zero dimensional case} \label{section 4}
As in the previous section, consider the weight function \WT\ on $k[x,y,z]$ defined by $\WT(x)=1$, $\WT(y)=2$, $\WT(z)=3$. We fix an integer $d \geq 5$ congruent to
$2$ modulo $3$, and we write $d=3m-1$ with $m\geq 2$.  Let
\[
g_0(x,y,z)= z^m+a_1(x,y)z^{m-1}+\cdots + a_{m-1} (x,y)z + a_m(x,y)
\]
denote a \WT-homogeneous polynomial in $k[x,y,z]$ of weight $3m$, monic in $z^m$.
\bp \label{prop-good-cohom}
Let $D$ be  the effective Cartier divisor on $L_{d-1}$ supported at the point $p_0=[0:0:0:1]$ with local equation $g_0$.
Consider as in Proposition \ref{liftDtoC} the line $L$ of equation $X=Y=0$ and the primitive $d$-structure $C_t$ on $L$
contained in the surface of equation $XGZ^a-tYW^{a+d-1}=0$. If $\HH^{0} (\ideal{D,\pso}(d-2))=0$,
then for general $t$ the sheaf $\caf(C_t) =\ideal{C_t}/\ideal{L}^d$ has natural cohomology, that is,
\[
\HH^0 \big(\ideal{C_t}(n)\big)= \HH^{0} \big(\ideal{L}^d(n)\big) \qquad \text{for all~}n \leq a+d-1.
\]
\ep

\begin{proof}
By Proposition \ref{hp}
\[
\chi \caf(n)=        \left(n-e-\frac{1}{3}(2d-1)\right)\binom{d}{2}
\]
hence $\chi \caf(n_0)=0$ for
\[
n_0= e+  \frac{1}{3}(2d-1)= a+d-1.
\]
Thus what we have to show is the vanishing $\HH^0 (\caf(a+d-1))=0$. By Lemma \ref{lemma} the sheaf
$\caf(a+d)$ is isomorphic to $\coo_{L_{d-1}}(E)$, and by semicontinuity, it is enough to prove $\HH^0 (\coo_{L_{d-1}}( E_{\infty}-H))=0$ where $E_{\infty}$ is the divisor
on $L_{d-1}$ supported at $p_{\infty}$ with local equation $g_{\infty}$. Next we note that the polynomial
$G$ defines a Cartier divisor $B$ on $L_{d-1}$ linearly equivalent to $(d-1)H$, where $H$ denotes the divisor
class of the plane section. The divisor $B$ is supported at $p_{\infty}$, where it has local equation $g_{\infty}$,
and at $p_0=[0:0:0:1]$ where it has local equation $g_0$. So let $D=B-E_{\infty}$ the divisor supported at $p_0$
with local equation $g_0$. Then
\[
\HH^0 \big(\coo_{L_{d-1}}( E_{\infty}-H)\big)= \HH^0 \big( \coo_{L_{d-1}} ((d-2)H-D)\big)=\HH^0 \big(\ideal{D,L_{d-1}} (d-2)\big) = 0. \qedhere
\]
\end{proof}

Finally, we note

\bl \label{grobner}
Let $D \subseteq \pso$ be a zero dimensional subscheme supported at  $p_0=[0:0:0:1]$ of length equal to $\HH^0 (\coo_{\pso}(d-2))$.
Suppose $I$ is the ideal of $D$ in $\mathbb{A}^3 = \pso \setminus \{W=0\} \cong \textnormal{Spec}( k[x,y,z])$. Then $D$ does not lie on a surface of degree $d-2$
if, and only if, $I$ contains no nonzero polynomial of degree $d-2$ or less. If we
fix a term order $\succcurlyeq$ on the monomials of $k[x,y,z]$ that refines the total degree partial order,
then $D$ does not lie in a surface of degree $d-2$ if and only if the initial ideal of $I$ with respect to $\succcurlyeq$
is $(x,y,z)^{d-1}$.
\el

\begin{proof}
Suppose that $S$ is a surface of minimal degree among those containing  $D$, and let $F=0$ be its equation, so that $F \in k[X,Y,Z,W]$ is a homogeneous polynomial, and we
let $s= \deg (F)$.
Since $D$ is supported at the origin of $\mathbb{A}^3$, no component of $D$ is contained in the plane at infinity $W=0$, so
that $F$ is not divisible by $W$ by minimality of its degree. Let $f=F(x,y,z,1)$: then $f=0$  is the affine equation of $S \cap \mathbb{A}^3$ in  $\mathbb{A}^3$,
and $f$ has degree $s$ because $F$ is not divisible by $W$. The initial term of $f$ is a form of degree $s$  because $\succcurlyeq$ refines
the total degree partial order. We conclude that $D$ does not lie in a surface of degree $d-2$ if and only if the initial ideal of $I$
does not contain a form of degree $s<d-1$; since the length of $D$ is
\[\HH^0 (\coo_{\pso}(d-2))= \dim k[x,y,z]/(x,y,z)^{d-1},\] the latter condition
is equivalent to the initial ideal being precisely $(x,y,z)^{d-1}$.
\end{proof}

The proof of the main Theorem \ref{main result} is now in place: statement \emph{(\ref{main:a})} is proven by Proposition \ref{liftDtoC}, statement \emph{(\ref{main:b})} follows Proposition \ref{prop-good-cohom}, Lemma \ref{grobner}.

\section{The polynomial $g$} \label{section gpol}
In this section we rephrase Lemma \ref{grobner}  in terms of the polynomial $g_0$ and we explain how one can check the validity of Conjecture \ref{conjB} for a given $m$ by computing the determinant of a (very large) matrix of linear forms. For $m = 2$ and $m = 3$, one can easily write such matrix for any polynomial $g$ of weight $3m$, while for bigger $m$ we consider a random polynomial and we verify Conjecture \ref{conjB} with the aid of \emph{Macaulay2} \cite{M2}.

We will use the following notation. A $\Z$-graded $k$-vector space $N$ is a vector space together with a decomposition $N= \bigoplus_{n \in \Z} N_n$
where each $N_n$ is a $k$-vector subspace. We will say the elements of $N_i$ are homogeneous of weight $n$. For every integer $q$, the $\Z$-graded $k$-vector space $N[q]$
it is just $N$ with weights shifted as follows $N[q]_n=N_{q+n}$.

Fixed an integer $m \geq 2$, we set $\ell=3m-2$ (so $\ell=d-1$ if $d=3m-1$ as in the previous section).  We denote by $R$ the $k$-algebra
$k[x,y,z]/(x,y)^{\ell}$ and,  by abuse of notation, we denote  the class of the monomial $x^i y^j z^k$ in $R$ still as $x^i y^j z^k$. A $k$-basis of $R$ consists of the monomial classes by $x^i y^j z^k$ for which
$i+j < \ell$. We will consider $R$ as a $\Z$-graded
vector space setting $R_n$ to denote the set of linear combinations of monomials of weight $n$,  where the weight function is defined by $\WT(x)=1$, $\WT(y)=2$ and $\WT(z)=3$.
We let $M$ denote the $\Z$-graded $k$-vector space obtained as the quotient of $R$ by the $\Z$-graded $k$-vector subspace generated
by monomials of standard degree $<\ell$.  We denote by $[h(x,y,z)]$ the image of the polynomial $h$ in $M$. Thus
a $k$-basis of $M$ is given by the classes
$[x^a y^b z^c]$  that satisfy $a+b<\ell$ and $a+b+c \geq \ell$. Note that such monomials have weight at least $3m$, and only one of them,
namely $[x^{\ell-1}z]$, has weight exactly $3m$. It turns out that $M$ can be given the structure of a weighted $R$-module that makes it
isomorphic to $R$ with generator  $[x^{\ell-1}z]$ if one defines $x [x^{\ell-1}z]=[x^{-1}y x^{\ell-1}z] $, $y [x^{\ell-1}z]=[x^{-1}z x^{\ell-1}z] $
and $z [x^{\ell-1}z]=[z x^{\ell-1}z] $. More precisely, we define a $k$-linear map
$\psi:R[-3m] \rightarrow M$  by the formula
\[
\psi(x^i y^j z^k)= [x^{\ell-i-j-1}y^{i}z^{1+j+k}]   \qquad (i+j < \ell).
\]
The map $\psi$ is well defined as the class in $R$ of the monomial $x^i y^j z^k$ is non zero if and only if $i+j < \ell$; it is \WT-homogeneous
of weight $3m$, that is,  $\WT (\psi (x^i y^j z^k))= \WT( x^i y^j z^k) + 3m$, and it is bijective, its inverse being
\[
\psi^{-1}([x^a y^b z^c])= x^{b}y^{\ell-1-a-b}z^{a+b+c -\ell}   \qquad (a+b<\ell \text{~and~}a+b+c \geq \ell).
\]
Thus $M \cong R[-3m]$ as $\Z$-graded $k$-vector spaces.

\bp \label{same-hilb}
Fix an integer $m \geq 2$ and define $\ell=3m-2$.
Let 
\[
g(x,y,z)=  z^m+a_1(x,y)z^{m-1}+\cdots + a_{m-1} (x,y)z + a_m(x,y)
\]
be a \WT-homogeneous polynomial in $k[x,y,z]$ of weight $3m$, monic in $z^m$. Consider in $k[x,y,z]$ the \WT-homogeneous
ideals $I=(x,y)^\ell+(g)$ and $J=(x,y,z)^\ell$. Then $k[x,y,z]/I$ and $k[x,y,z]/J$ are isomorphic as $\Z$-graded $k$-vector
spaces.
\ep
\begin{proof}
Let $S=k[x,y,z]/I$ and $T=k[x,y,z]/J$. We need to show $\dim S_n=\dim T_n$ for every $n$. For this is enough to show
there are monomial bases of $S$ and $T$ with the same number of elements of weight $n$ for every $n$. Now a monomial basis
$\mathcal{B}(S)$ of $S$ is  given by the classes of the monomials $x^a y^b z^c$ that satisfy $a+b < \ell$ and $c<m$, while
a monomial basis $\mathcal{B}(T)$ of $T$ is  given by the classes of the monomials $x^i y^j z^k$  that satisfy $i+j+k < \ell$.

Consider the map $\phi_{z^m}:R[-3m] \rightarrow M$ defined by multiplication by $z^m$ followed by the $k$-linear projection from $R$
to $M$ that kills all monomials of standard degree $<\ell$:
\[\phi_{z^m}(h)=[z^m h].\]
Let $K$ and $Q$ the kernel and the cokernel of $\phi_{z^m}$: these are $\Z$-graded $k$-vector spaces, and as such they are isomorphic
because  $R[-3m] \cong M$ and $\phi_{z^m}$ is homogeneous. A monomial basis of the cokernel $Q$ consists of those monomials $[x^a y^b z^c]$
that satisfy $c <m$ besides $a+b<\ell$ and $a+b+c \geq \ell$; these are the monomials of the basis $\mathcal{B}(S)$ that are not
in $\mathcal{B}(T)$. On the other hand, a monomial basis of the kernel $K$, consists of those monomials $x^i y^j z^k$ that satisfy $i+j+k+m < \ell$;
multiplying these monomials by $z^m$ we obtain the monomials of the basis $\mathcal{B}(T)$ that are not
in $\mathcal{B}(S)$. Hence $\mathcal{B}(S)$ has the same number of elements of weight $n$ as  $\mathcal{B}(T)$, for every $n$.
\end{proof}

In the following proposition we sum up the equivalent conditions we need on the polynomial $g$ for the curves constructed  in Section  \ref{section 4} with $g_0=g$
to have good cohomology.
\bp \label{good g}
Fix an integer $m \geq 2$ and let $\ell=3m-2$.
Let $g$ be a \WT-homogeneous polynomial in $k[x,y,z]$ of weight $3m$. The following conditions are equivalent:
\ben[(i)]
\item\label{good g:1}
The ideal $I=(x,y)^\ell+(g)$ contains no polynomial of standard degree $\ell-1$.
\item\label{good g:2}
If  $\succcurlyeq$ is a term order on the monomials of $k[x,y,z]$ that refines the standard degree partial order,
the initial ideal of $I=(x,y)^\ell+(g)$ with respect to $\succcurlyeq$ is $(x,y,z)^{\ell}$.
\item\label{good g:3}
 The \WT-homogeneous $k$-linear map $\phi_{g}:R[-3m] \rightarrow M$ defined by multiplication by $g$ followed by the $k$-linear projection from $R$
to $M$: \[\phi_{g}(h)=[g h]\]
is an isomorphism.
\een
\ep
\begin{proof}
By Lemma \ref{grobner}, the first two conditions are equivalent. It is clear that $\phi_g$ is an isomorphism if, and only if, the ideal $I=(x,y)^\ell+(g)$ contains no
polynomial of standard degree $\ell-1$.
\end{proof}

\br\label{finite cases}
Assume $z^m$ appears in $g$ with nonzero coefficient. In order to show $\phi_{g}:R[-3m] \rightarrow M$ is an isomorphism, it is enough to show that
\[
(\phi_{g})_{n}: R[-3m]_{n} \rightarrow M_{n}
\]
is an isomorphism for $3m \leq n \leq 9(m-1)$. Indeed,  every monomial of weight $n=a+2b+3c > 9(m-1)$ with $a+b < \ell=3m-2$ satisfies $c \geq m$,
hence is in the image of $\phi_{z^m}$. In other terms, writing as above $J=(x,y,z)^\ell$ and $I=(x,y)^\ell+(g)$, the isomorphic $\Z$-graded vector spaces
$k[x,y,z]/I$ and $k[x,y,z]/J$ have no element weight $>9(m-1)$. In particular, every polynomial of weight $>9(m-1)$ is a multiple of $g$ modulo $(x,y)^\ell$,
therefore  $(\phi_{g})_{n}$ is surjective, hence an isomorphism, for every $n>9(m-1)$.
\er

\bex \label{m23}
Consider $m=2$ and let $\ell=3m-2 = 4$. A polynomial $g$ of weight $3m = 6$ in $k[x,y,z]$ has the form
\[
g(x,y,z) = a_0 z^{2}+ a_1  x^{3} z+ a_2 x y z +a_{3} y^{3}  \qquad (a_0, a_1, a_2, a_3 \in k).
\]
We want to prove Conjecture \ref{conjB} by describing explicitly for which values $a_i$, the ideal $I=(x,y)^\ell+(g)$ contains no polynomial of standard degree $\ell-1$. Using the equivalent condition \emph{(\ref{good g:3})} of Theorem \ref{good g} and taking into account Remark \ref{finite cases}, we consider the morphism $(\phi_{g})_{n}$ for $6 \leq n \leq 9$:
\[
\begin{split}
&(\phi_{g})_{6} : R[-6]_6 = \langle 1 \rangle \xrightarrow{\scriptsize\left[\begin{array}{c} a_1 \end{array}\right]} M_6 = \langle x^3 z\rangle,\\
&(\phi_{g})_{7} : R[-6]_7 = \langle x \rangle \xrightarrow{\scriptsize\left[\begin{array}{c} a_2 \end{array}\right]} M_7 = \langle x^2 y z\rangle,\\
&(\phi_{g})_{8} : R[-6]_8 = \langle x^2, y \rangle \xrightarrow{\scriptsize\left[\begin{array}{cc} a_0 & 0 \\ 0 & a_2 \end{array}\right]} M_8 = \langle x^2 z^2, x y^2 z\rangle,\\
&(\phi_{g})_{9} : R[-6]_9 = \langle z, x^3, x y \rangle \xrightarrow{\scriptsize\left[\begin{array}{ccc} a_1 & a_2 & a_3\\ a_0 & 0 & 0 \\ 0 & a_0 & 0 \end{array}\right]} M_9 = \langle z^2 x^3 , z^2 x y, z y^3 \rangle.\\
\end{split}
\]
Therefore  $\phi_{g}$ is an isomorphism provided that all coefficients of $g$ are non zero.

The calculation becomes much longer for $m=3$. In this case the polynomial $g$ has weight $9$ and the form
\[
  g(x,y,z)= {b}_{0} z^{3}+{b}_{1} x^{3} z^{2} +{b}_{2} x y z^{2} +{b}_{3} x^{6} z +{b}_{4} x^{4} y z +{b}_{5}
      x^{2} y^{2} z +{b}_{6} y^{3} z +{b}_{7} x^{3} y^{3}+{b}_{8} x y^{4}.
\]
One needs to check that $(\phi_{g})_{n}$ is an isomorphism for $9 \leq n \leq 18$.  We spare the reader the details (but the interested reader can find all details in the ancillary \emph{Macaulay2} file \texttt{MaximumGenusProblem.m2} available at the webpage \url{www.paololella.it/EN/Publications.html}). In this case, $b_0$, $b_1$, $b_3$, $b_4$ and $b_6$ are required to be non zero, and there are other 5 non-vanishing conditions of degree 2, 3, and 4 among the coefficients $b_i$ describing $g$. For instance, the condition $(\phi_{g})_{17}$ is an isomorphism is $b_0^6(2b_2 b_6- b_0 b_8) \neq 0$. Similarly for every $9 \leq n \leq 18$ the set of $g$ for which  $(\phi_{g})_{n}$ is not an isomorphism is a hypersurface in $\mathbb{A}^8$. Thus $\phi_{g}$ is an isomorphism if $g$ is general.
\eex

For larger $m$, it is more convenient to prove Conjecture \ref{conjB} directly, i.e.~by exhibiting a polynomial $g$ satisfying condition \emph{(\ref{good g:1})} of Theorem \ref{good g}. Indeed, if a special polynomial satisfies Theorem \ref{good g}\emph{(\ref{good g:1})}, then also the general $g$ does. We use condition \emph{(\ref{good g:2})} of Theorem \ref{good g} and we proceed as follows. With the aid of \emph{Macaulay2}, we construct a \WT-homogeneous polynomial $g$ with randomly chosen coefficients and we compute the initial ideal of $(x,y)^{\ell}+ (g)$ with respect to the graded lexicographic ordering. In order to speed up the computation, we consider polynomials with coefficients in a finite field $\mathbb{Z}/p\mathbb{Z}$, so that conjectures turn out to be proved for fields with sufficiently large characteristic (and surely for characteristic 0). In a reasonable computing time, it is possible to verify Conjecture \ref{conjB} for $m \leq 40$, that is Conjecture \ref{A} for $d \leq 120$. The interested reader can find the code in the aforementioned \emph{Macaulay2} file \texttt{MaximumGenusProblem.m2}.

%
%
%

\section{Sharpness of $P(d,s)$} \label{section 5}
In this section we show how to deduce sharpness of the bound $P(d,s)$ 
for $d \geq 2s-1$ if one knows $P(s{-}1,s{-}1)$ is sharp. We also briefly speculate on the intermediate
range $s+1 \leq d \leq 2s-2$.

\bp If $P(s{-}1,s{-}1)= P_{\textsc{max}} (s{-}1,s{-}1)$, then 
$P(d,s)=  P_{\textsc{max}}(d,s)$ for every  $d \geq 2s-1$.
\ep

\begin{proof}
We prove the statement by biliaison: since $P(s{-}1,s{-}1)$ is sharp, we can find
a curve $Y$ of degree $s\!-\!1$, not lying
on a surface of degree $s\!-\!2$, of maximum genus $g(Y)=P(s{-}1,s{-}1)$. If
$d \geq 2s\!-\!1$ and we set $t=d-s+1$, then $ t \geq s$ and therefore
the curve $Y$ lies in a surface $S$ of degree $t$.
Let $C$ any curve linearly equivalent to $Y+H$ on $S$. Then
the genus of $C$ is \cite[p.~66]{MDP}
\[
g(C)=g(Y)+1/2(t-1)(t-2)+\deg(Y)-1
\]
which coincides with $P(d,s)$
provided $d \geq 2s-1$; furthermore, $C$ does not lie on a surface of degree $s\!-\!1$ because $t\geq s$.
\end{proof}

\br 
To study sharpness of $P(d,s)$ in the range $s+1 \leq 2s-2$ one needs something different. The following remark might help.
Recall that, if a curve $C$ is the union of two curves $X$ and $Y$ that meet in a zero dimensional scheme $X\centerdot Y$ of length $\ell$,
then
\[
p_a (C)= p_a (X) + p_a (Y) + \ell -1.
\]
Therefore:
\begin{enumerate}[a)]
      \item suppose $s<d \leq 2s$ and $D$ is a curve of degree $d-1$ and arithmetic genus $p_a(D)=P(d-1,s)$. Then, if
  $L$ is a line that meets $D$ in a scheme of length $s$, then
  the curve $C=D \cup L$ has degree $d$ and arithmetic genus $P(d,s)$.
    \item\label{b} Suppose $s<d \leq 2s$ and $Y$ is a curve of degree $s$ and arithmetic genus $p_a(Y)=P(s,s)$. Then, if
  $X$ is a plane curve of degree $k=d-s$ that meets $Y$ in a scheme of length
  \[\ell= s+(s-1) + \cdots +(s-k+1)= ks-\frac{1}{2}k(k-1),\]
  the curve $C=X \cup Y$ has degree $d$ and arithmetic genus $P(d,s)$.
    \item\label{c} Suppose $d=s+k$ with $3 \leq k \leq s-2$. Let $Y$ denote a curve of degree $s-1$ and arithmetic genus $p_a(Y)=P(s{-}1,s{-}1)$,
  and let $P$ be a plane curve of degree $k+1$. Then  $C=P \cup Y$ has arithmetic genus $p_a(C)=P(d,s)$ provided
  \[
  P\centerdot Y= (s-1)-\frac{1}{2} (k-s)(k-s+1).
  \]
\end{enumerate}
The trouble in \ref{b}) is that we don't know how to construct such an  $X$ if $k \geq 3$, while \ref{c}) can be useful only for $k$ close to $s$.
\er

As an application we can prove
\bc
Bound \eqref{Bbound} is sharp for $s=5$, namely, $P(d,5)=  P_{\textsc{max}}(d,5)$ holds for every  $d \geq 5$.
\ec

\begin{proof}
By our previous results $P(d,5)=  P_{\textsc{max}}(d,5)$ if $d=5$ or if $d \geq 9$. It remains to check $d=6,7,8$.
We have seen how to construct a primitive $5$ line $C$ of genus $P(5,5)=-14$ that does not lie on a surface of degree $4$;
$C$ lies on a quintic surface that has an equation of the form $XG-YW^4=0$. The surface $S$ contains, besides the line that supports $C$, 
the two lines $L'$ and $L''$ of equations respectively $Y=Z=0$ and $X=W=0$. The planes $Y=0$ and $X=0$ are tangent to $C$. We let
$D'=V(Y,Z^2)$ and $D''=V(X,W^2)$. Then  one can check that $L'$ and $L''$ cut $C$ in a scheme of length $5$, while $D'$ and $D''$   
cut $C$ in a scheme of length $9$. Hence the curves $C_6= C \cup L'$,  $C_7= C \cup D'$, $C_8= C \cup D' \cup L''$ 
satisfy $p_a(C_d)=P(d,5)$ for $d=6,7,8$. 
\end{proof}

%
%

\end{document}